\documentclass[12t]{amsart}

\usepackage[pagebackref,colorlinks=true]{hyperref}  

\usepackage{eucal,url,amssymb,verbatim, 
enumerate,amscd,
}

\numberwithin{equation}{section}

\newtheorem{thrm}{Theorem}[section]

\newtheorem{lemma}[thrm]{Lemma}

\newtheorem{cor}[thrm]{Corollary}

\newcommand{\R}{\mathbb{R}}

\def\sideremark#1{\ifvmode\leavevmode\fi\vadjust{\vbox to0pt{\vss
 \hbox to 0pt{\hskip\hsize\hskip1em
 \vbox{\hsize2.5cm\tiny\raggedright\pretolerance10000
 \noindent #1\hfill}\hss}\vbox to8pt{\vfil}\vss}}}%

\renewcommand{\sideremark}[1]{}

\begin{document}

\begin{abstract}
We show that the CR structure on the twistor space of a quaternionic contact structure described by
O. Biquard is normal if and only if the Ricci curvature of the Biquard connection commutes with the
endomorphisms in the quaternionic structure of the contact distribution.
\end{abstract}



\title[The twistor space of a quaternionic contact manifold]
{The twistor space of a quaternionic contact manifold}
\date{\today}

\author{Johann Davidov}
\address[Johann Davidov]{ Institute of Mathematics and Informatics, Bulgarian Academy of Sciences.
1113 Sofia, Bulgaria \\
and "L. Karavelov" Civil Engineering Higher School,1373 Sofia, Bulgaria} \email{jtd@math.bas.bg}

\author{Stefan Ivanov}
\address[Stefan Ivanov]{University of Sofia, Faculty of Mathematics and Informatics,
blvd. James Bourchier 5, 1164,
Sofia, Bulgaria} \email{ivanovsp@fmi.uni-sofia.bg}

\author{Ivan Minchev}
\address[Ivan Minchev]{University of Sofia\\
Sofia, Bulgaria\\
and  Mathematik und Informatik\\
Philipps-Universit\"at Marburg\\
Hans-Meerwein-Str. / Campus Lahnberge 35032 Marburg, Germany}
\email{minchevim@yahoo.com}

\maketitle
\tableofcontents


\setcounter{tocdepth}{2}

\section{Introduction}
The notion of a  quaternionic contact (QC) structure is introduced in \cite{Biq1} and it describes
a type of geometrical structure that appears naturally as the conformal boundary at infinity of the
quaternionic hyperbolic space. {In general, a  QC structure on a real (4n+3)-dimensional manifold
$M$ is a codimension
three distribution $H$, the contact distribution , locally given as the kernel of { a} 1-form $%
\eta=(\eta_1,\eta_2,\eta_3)$ with values in $\mathbb{R}^3$ such that the three 2-forms $d\eta_i|_H$
are the fundamental 2-forms of a quaternionic structure on $H$ (for more details see the next
section).

It is a fundamental theorem of Biquard \cite{Biq1} that a QC-structure on a real analytic manifold
$M^{4n+3}$ is always the conformal infinity of a quaternionic K\"ahler metric defined in the
neighborhood of $M^{4n+3}$. This theorem generalizes an earlier result of LeBrun  \cite{LeB} that
states that a real analytic conformal 3-manifold is always the conformal infinity of a self-dual
Einstein metric. From this point of view, one may regard the QC-geometry as a natural
generalization to  dimensions $4n+3$ of the 3-dimensional conformal Riemannian geometry. Moreover,
the QC-geometry gives a natural setting for certain Yamabe-type problem \cite{Wei,IMV1,IMV2,IMV4}.
A particular case of this problem amounts to finding the extremals and the best constant in the
$L^2$ Folland-Stein Sobolev-type embedding, \cite{F2} and \cite{FS}, on the quaternionic Heisenberg
group, see \cite{GV} and \cite{IMV2, IMV4}.

The 1-form $\eta$ that defines the QC-structure is determined up to a conformal factor and the
action of $SO(3)$ on $\mathbb{R}^3$. {Therefore} $H$ is equipped with a conformal class $[g]$ of
metrics and a 3-dimensional  quaternionic bundle $Q$. The associated 2-sphere bundle
$S^2(Q)\rightarrow M$ is called the twistor space of the QC-structure. The transformations
preserving given QC structure $\eta$, i.e. the transformations of the type
$\bar\eta=\mu\Psi\cdot\eta$ for a positive smooth function $\mu$ and an $SO(3)$ matrix $\Psi$ with
smooth functions as entries, are called \emph{quaternionic contact conformal (QC conformal)
transformations}. If the function $\mu$ is constant we have {\em quaternionic contact homothetic
(QC homothetic) transformations}. To every metric in the fixed conformal class $[g]$ on $H$ one can
associate a linear connection preserving the QC structure, \cite{Biq1}, which we shall call the
Biquard connection. This connection is invariant under QC homothetic transformations but changes in
a non-trivial way under QC conformal transformations.

Examples of QC manifolds can be found in \cite{Biq1,Biq2,IMV1,D1}. It is known that on the sphere
$S^{4n+3},\ n>1$  there exist infinitely  many different global QC-structures. Indeed, in
\cite{LeB2} Lebrun has constructed an infinite-dimensional space of deformations of the standard
hyperbolic quaternionic-K\"ahler metric on the open Ball $B^{4n+4}$ through complete
quaternionic-K\"ahler metrics. After this,  Biquard \cite{Biq1} has shown that each of the
constructed metrics actually has as a conformal infinity a certain unique QC-structure on the
boundary $S^{4n+3}$ of the ball. However his construction does not give the QC-structures
explicitly. The amount of known explicit examples remains very restricted.

As we have mentioned  above and  shall explain  in the next paragraph, each QC-structures with a
fixed metric in the conformal class $[g]$ determines a unique connection $\nabla,$  the Biquard
connection. This connection plays a roll in the QC-geometry similar to the one played by the
Levi-Civita connection in the 3-dimensional conformal geometry. The restriction to $H$ of the Ricci
tensor of $(g,\nabla)$ gives rise to three quantities, namely the QC-scalar curvature $Scal$ and
two symmetric trace-free (0,2) tensor fields $T^0$ and $U$ defined on the contact distribution $H$.
The tensors $T^0$ and $U$ determine the trace-free part of the Ricci tensor restricted to $H$ and
can also be expressed in terms of the torsion endomorphisms of the Biquard connection \cite{IMV1}
(see Section 2 for the details). According to \cite{IMV1}, the vanishing of the torsion
endomorphisms of the Biquard connection is equivalent to $T^0=U=0$ and if the dimension is at least
eleven, then the function $Scal$ has to be constant. If in addition this constant is different from
zero, then the QC-structure is locally QC-homothetic to a (positive or negative) $3$-Sasakian
structure (see also \cite{IV1,IMV3}).

{ Explicit examples of QC manifolds with zero or non-zero torsion endomorphism  have been recently
given in \cite{AFIV,AFIV1}}. The quaternionic Heisenberg group, the quaternionic sphere of
dimension $4n+3$ with its standard 3-Sasakian structure and the QC structures locally QC conformal
to them are characterized in \cite{IV} by the vanishing of a tensor invariant under conformal
transformations, the QC-conformal curvature defined in terms of the curvature and the torsion of
the Biquard connection. Explicit examples of non QC conformally flat QC manifolds are constructed
in \cite{AFIV,AFIV1}.

The twistor space ${\mathcal Z}=S^2(Q)$ of a QC-structure is naturally equipped with a CR-structure
\cite{Biq1,D}, which is invariant under the QC-conformal transformations. To each metric $g\in[g]$
one can naturally define a contact form $\eta^{\mathcal Z}$ on the twistor space ${\mathcal Z}$
compatible with the CR-structure there. The contact form depends on the choice of $g\in[g]$ and
thus the whole construction is not  QC-conformal invariant anymore but it remains however
QC-homothetic invariant.

The purpose of the present notes is to show (Theorem~\ref{main_thrm}) that the contact form
$\eta^{\mathcal Z}$ is normal if and only if the tensor $T^0$ vanishes. The latter condition is
equivalent to the requirement that the Ricci tensor of the Biquard connection commutes with the
endomorphisms in the quaternionic structure of $H$. Note that the normality of the contact manifold
$({\mathcal Z},\eta^{\mathcal Z})$ is equivalent to the condition that the product manifold
${\mathcal Z}\times \R$ is a complex manifold with a certain naturally defined complex structure
(cf. e.g. \cite{Bl}).

Note also that every (negative or positive) 3-Sasakian manifold has constant QC-scalar curvature
and satisfies the condition $T^0=U=0$ \cite{IMV1}. According to Theorem~\ref{main_thrm} the
CR-structure on its twistor space is normal.
 It is shown in \cite{IMV3} that, in the  case of zero torsion endomorphisms of the Biquard connection,
 the vector bundle $Q\rightarrow M$ admits a flat connection which implies that the corresponding bundle $\tilde Q\rightarrow \tilde M$
of the universal cover $\tilde M$ of $M$ is trivial. Thus in the case $T^0=U=0$ (plus the condition
$Scal$=const in dimension seven) we see that the twistor space $\tilde{\mathcal Z}$ of the
universal cover of $M$ is just the product $\tilde{\mathcal Z}=\tilde M\times S^2.$

To the best of our knowledge no explicit examples of QC structures with  $T^0=0$  and non-constant
QC scalar curvature are known.

\smallskip

\textbf{Acknowledgements} S.I. and I.M. are partially supported by the Contract 181/2011 with the
University of Sofia `St.Kl.Ohridski', Contracts ``Idei", DO 02-257/18.12.2008 and DID
02-39/21.12.2009. J.D is partially supported by "L.Karavelov" Civil Engineering Higher School,
Sofia, Bulgaria under contract No 10/2009.

\section{Quaternionic contact manifolds and the Biquard connection}
\label{s:review} In this section we will briefly review the basic notions of quaternionic contact
geometry and recall certain results of \cite{Biq1} and \cite{IMV1}.

A quaternionic contact structure (shortly, QC-structure) on a (4n+3)-dimensional smooth manifold
$M$ consists of a rank $4n$ subbundle $H$ of $TM$, a positive definite metric $g$ on $H$ and a rank
3 subbundle $Q$ of $End(H)$ such that, in a neighbourhood $U$ of each point of $M$, there are
$1$-form $\eta=(\eta_1,\eta_2,\eta_3)$ with values in $\R^3$ and a triple $\vartheta=(I_1,I_2,I_3)$
of sections of $Q$ with the following properties:


\par (1) $H|U$ is the kernel of $\eta$;


\par (2) The bundle $Q$ is locally generated
by three almost complex structures $I_1,I_2,I_3$ on $H$ satisfying the identities of the imaginary
unit quaternions, $I_1^2=I_2^2=I_3^2=-Id_{H}$, $I_1I_2=-I_2I_1=I_3$.



\par (3)  $d\eta_s(X,Y)=2g(I_sX,Y)$ for $X,Y\in H|U$.


\medskip

{\bf Convention}.

\begin{enumerate}
\item[a)] Throughout this paper, we shall use $X,Y,Z,U$ to denote vectors or sections of $H$;
\item[b)] $\{e_1,\dots,e_{4n}\}$ denotes a  local orthonormal basis of $H$;
\item[c)] The triple $(i,j,k)$ denotes any cyclic permutation of $(1,2,3)$.
\item[d)] $s,t$  will be any number from the set $\{1,2,3\}, \quad
s,t\in\{1,2,3\}$.
\end{enumerate}

If $\eta=(\eta_1,\eta_2,\eta_3)$ and $\vartheta=(I_1,I_2,I_3)$ satisfy conditions (1), (2), (3), we
shall say that $(\eta,\vartheta)$ is an {\it admissible set} for the $QC$-structure.


Any two triples of sections of $Q$ satisfying condition (2) constitute frames of $Q$ which induce
the same orientation, thus the bundle $Q$ has a canonical orientation.

Condition $(3)$ implies that $g(IX,Y)=-g(X,IY)$ for any section $I$ of $Q$ and $X,Y\in H$ and,
thus, $H$ is equipped with an $Sp(n)Sp(1)$-structure.


The metric $g$ induces a metric on the bundle $End(H)$ defined by
$
<A,B>=\frac{1}{4n} Trace A^t B,
$
where $A$ and $B$ are endomorphisms of a fibre of $H$ and $A^t$ is the adjoint of $A$ with respect
to $g$. Any sections  $I_1,I_2,I_3$ of  $End(H)$ satisfying the imaginary quaternion relations form
an orthonormal set with respect to the induced metric.  Moreover, if $I$ and $J$ are sections of
$Q$, then $<I,J>=0$ if and only if $IJ=-JI$ and $<I,I>=1$ exactly when $I^2=-Id_H$.  Note also that
any oriented orthonormal frame of $Q$ consists of endomorphisms of $H$ satisfying the imaginary
quaternion relations.


Given a distribution $H$ on a smooth manifold $M$ and a vector bundle $E$ over $M$, a partial
connection on $E$ along $H$ is, by definition, a bilinear map $\nabla_{X}\sigma$ defined for vector
fields $X$ with values in $H$ and sections $\sigma$ of $E$ such that
$\nabla_{fX}\sigma=f\nabla_{X}\sigma$ and $\nabla_{X}(f\sigma)=X(f)\sigma+f\nabla_{X}\sigma$ for
every smooth function $f$ on $M$.

Let $H$ be a distribution of a manifold $M$ and $g$ be a metric on $H$. Biquard \cite[Lemma
II.1.1]{Biq1} has observed that, for any supplementary distribution $V$ of $H$ in $TM$,  there is a
unique partial connection $\nabla$ on $H$ along $H$ such that
\begin{enumerate}
\item[$(i)$] $\nabla g = 0$;
\item[$(ii)$] for any two sections $X,Y$ of $H$, the torsion $T(X,Y)=\nabla_{X}Y-\nabla_{Y}X-[X,Y]$ satisfies the identity
$T(X,Y)=-[X,Y]_V$, where the subscript $V$ means "the component in $V$";
\end{enumerate}

  Now let $(M,H,g,Q)$ be a quaternionic contact manifold. Fix a supplementary distribution $V$ of $H$ in
$TM$ and let $\nabla$ be the associated connection on $H$ along $H$. The partial connection on
$End(H)$ along $H$ induced by $\nabla$ will be denoted also by $\nabla$. Biquard \cite[Lemma
II.1.6, Proposition II.1.7]{Biq1} has shown that $\nabla$ preserves the bundle $Q$ if and only if,
around any point of $M$, there is an admissible set $(\eta,\vartheta)$ such that the frame
$(\xi_1,\xi_2,\xi_3)$ of $V$ dual to the frame $(\eta_1|V, \eta_2|V, \eta_3|V)$ satisfies the
condition
\begin{equation}\label{bi1}
(\imath_{\xi_s}d\eta_{t})|H=-(\imath_{\xi_t}d\eta_{s})|H,\quad s,t=1,2,3.
\end{equation}
where $\imath$ denotes the interior multiplication. Note that, if  condition \eqref{bi1} is
satisfied for an admissible set $(\eta,\vartheta)$, then it holds for any other admissible set
$(\eta',\vartheta')$. Indeed, we have $\eta_t'=\sum_{s=1}^3 a_{ts}\eta_s$,  where $[a_{ts}]$ is a
non-singular $3\times 3$-matrix of smooth functions. In view of $(1)$,  $d\eta_t'(X,Y)=\sum_{s=1}^3
a_{ts}d\eta_s(X,Y)$ for $X,Y\in H$, hence $I_t'=\sum_{s=1}^3 a_{ts}I_s$ by $(3)$. The latter
identity and $(2)$ imply that $[a_{ts}]\in SO(3)$. Then $\xi_t'=\sum_{s=1}^3 a_{ts}\xi_s$, where
$(\xi_1',\xi_2',\xi_3')$ is the dual frame of $(\eta_1'|V,\eta_2'|V,\eta_3'|V)$. This observation
implies our claim.


Buqiuard \cite[Th\'eor\`{e}m II.1.3]{Biq1} has proved that if $dim\,M>7$, then their is a unique
supplementary distribution $V$ of $H$ in $TM$ for which the associated connection $\nabla$
preserves the bundle $Q$:

\begin{enumerate}
\item[$(iii)$]  $\nabla_{X}Q\subset Q$ for $X\in H$.
\end{enumerate}

For any admissible set $(\eta,\vartheta)$, the frame $(\xi_1,\xi_2,\xi_3)$ of the bundle $V$ dual
to the frame $(\eta_1|V, \eta_2|V, \eta_3|V)$ will be called {\it associated} to
$(\eta,\vartheta)$.

\smallskip

Given a section $\xi$ of $V$ and a section $X$ of $H$, set

\begin{enumerate}
\item[$(iv)$]  $\nabla_X\xi=[X,\xi]_V$.
\end{enumerate}

By \cite[Proposition II.1.9]{Biq1}, the latter formula defines a partial connection on $V$ along
$H$ such that
\begin{enumerate}
\item[$(v')$]  $\nabla <.,.>=0$
\end{enumerate}


Let $(\eta,\vartheta)$ be an admissible set for the given quaternionic contact structure on $M$ and
let $(\xi_1,\xi_2,\xi_3)$ be the frame of $V$ associated to $(\eta,\vartheta)$. Then the assignment
\begin{equation}\label{phi}
\xi_s\to I_s, s=1,2,3,
\end{equation}
determines an bundle isomorphism $\varphi:V\to Q$ that does not depend on the particular choice of
the admissible set. The isomorphism $\varphi$ has the property that $\nabla_{X}\varphi=0$ for $X\in
H$. Indeed, by $(3)$ and (\ref{bi1}), we have
$$
\begin{array}{c}
\nabla_X\varphi(\xi_t)=\nabla_X I_t=-\sum_{s=1}^3 d\eta_t(\xi_s,X)I_s=\sum_{s=1}^3
d\eta_s(\xi_t,X)I_s=\sum_{s=1}^3\eta_s([X,\xi_t]_V)I_s=\\[6pt]
\sum_{s=1}^3\eta_s(\nabla_X\xi_t)\varphi(\xi_s)=\varphi(\nabla_X\xi_t)
\end{array}
$$


Set
$$P=\{A\in End(H)\ |\ A ~\mbox{is skew-symmetric and}~AI=IA ~\mbox{for every}~I\in Q\}.$$
This is a subbundle of $End(H)$ of rank $2n^2+n$, orthogonal to $Q$ and such that the commutator
$[A_1,A_2]$ of two endomorphisms $A_1, A_2\in P$ is also in $P$. Clearly, every fibre of $P$ (resp.
$Q$) is isomorphic to the Lie algebra $sp(n)$ (resp. $sp(1)$).

 It is shown in \cite[Lemme II.2.1]{Biq1} that there is a unique partial connection $\nabla$ on $H$
along $V$ such that

\begin{enumerate}
\item[$(v)$] $\nabla g = 0$;
\item[$(vi)$] The induced connection on $End(H)$ preserves the bundle $Q$;
\item[$(vii)$] Setting $T(\xi,X)=\nabla_{\xi}X-\nabla_{X}\xi-[\xi,X]$ for $\xi\in V$ and $X\in H$,
every endomorphism $$T_{\xi}:\,H\in X\to T(\xi,X)=\nabla_{\xi}X-[\xi,X]_{H}\in H$$ is an element of
$(P\oplus Q)^{\bot}\subset End(H)$.
\end{enumerate}

Note, that we have a bundle isomorphism $\{(P\oplus Q)^{\bot}\subset End(H)\}\cong \{(sp(n)\oplus
sp(1))^{\bot}\subset gl(4n)\}$.


 Since $\nabla_{\xi}Q\subset Q$ for every $\xi\in V$, we can transfer $\nabla_{\xi}$ from $Q$
to $V$ via the isomorphism $\varphi:V\to Q$. In this way get a partial connection on $V$ along $V$.

Combining the partial connections we have defined, we obtain a connection $\nabla$ on $TM$ having
the properties $(i)$-$(vii)$ and the property

\noindent $(viii)$ $\nabla\varphi =0$.


We shall call $\nabla$ \emph{the Biquard connection} of the $QC$-structure $(H,g,Q)$ on $M$.


In the case when the dimension of $M$ is seven, it is not always possible to find a supplement $V$
to $H$ for which condition \eqref{bi1} hold. Duchemin \cite{D} has shown that if we assume that,
around any point of $M$, there exists an admissible set for which we can find vector fields
$\xi_1,\xi_2,\xi_3$ satisfying \eqref{bi1}, then one can define a connection with the properties
$(i)$-$(viii)$. Henceforth, by a quaternionic contact structure in dimension $7$ we shall mean a
QC-structure satisfying \eqref{bi1}.


Let $\varphi:V\to Q$ be the isomorphism defined by (\ref{phi}). Using this isomorphism we transfer
to $V$ the metric and the orientation of $Q$. Then any frame $\xi_1,\xi_2,\xi_3$ associated to an
admissible set of the $QC$-structure is orthonormal and positively oriented. Putting together the
metric of $V$ and the metric $g$ of $H$ we obtain a metric on $TM=H\oplus V$  for which $H$ and $V$
are orthogonal. This metric will be also denoted by $g$. It follows from properties $(v')$, $(v)$
and $(vi)$ that the connection $\nabla$ on $Q$ is compatible with the metric $<.,.>$. Therefore the
metric $g$ on $TM$ is parallel with respect to the Biquard connection, $\nabla g=0$.


The properties of the Biquard connection are encoded in the properties of the torsion endomorphisms
$T_{\xi}=T(\xi,\cdot) : H\rightarrow H, \quad \xi\in V$.

Any endomorphism $\Psi$ of $H$ can uniquely be  decomposed with respect to the quaternionic
structure $(Q,g)$ into four $Sp(n)$-invariant parts
$\Psi=\Psi^{+++}+\Psi^{+--}+\Psi^{-+-}+\Psi^{--+},$
where $\Psi^{+++}$ commutes with all three $I_i$, $\Psi^{+--}$ commutes with $I_1$ and
anti-commutes with the others two and so on. Explicitly,
\begin{equation*}
\begin{aligned}4\Psi^{+++}=\Psi-I_1\Psi I_1-I_2\Psi I_2-I_3\Psi
I_3,\quad
4\Psi^{+--}=\Psi-I_1\Psi I_1+I_2\Psi I_2+I_3\Psi I_3,\\
4\Psi^{-+-}=\Psi+I_1\Psi I_1-I_2\Psi I_2+I_3\Psi I_3,\quad 4\Psi^{--+}=\Psi+I_1\Psi
I_1+I_2\Psi I_2-I_3\Psi I_3.
\end{aligned}
\end{equation*}
\noindent The two $Sp(n)Sp(1)$-invariant components are $\Psi^{+++}$ and
$\Psi^{+--}+\Psi^{-+-}+\Psi^{--+}$. If $n=1$, then the space of symmetric endomorphisms commuting
with all $I_s$ is 1-dimensional, i.e. $\Psi^{+++}$  is proportional to the identity,
$\Psi^{+++}=\frac{|\Psi|^2}{4}Id_{|H}$.

Decomposing the endomorphism $T_{\xi}\in(sp(n)+sp(1))^{\perp}$ into its symmetric part $T^0_{\xi}$
and skew-symmetric part $b_{\xi}, T_{\xi}=T^0_{\xi} + b_{\xi} $, Biquard has shown in \cite{Biq1}
that the torsion $T_{\xi}$ is completely trace-free, $tr\, T_{\xi}=tr\, T_{\xi}\circ I_s=0$, its
symmetric part has the properties $T^0_{\xi_i}I_i=-I_iT^0_{\xi_i}\quad
I_2(T^0_{\xi_2})^{+--}=I_1(T^0_{\xi_1})^{-+-},\quad
I_3(T^0_{\xi_3})^{-+-}=I_2(T^0_{\xi_2})^{--+},\quad I_1(T^0_{\xi_1})^{--+}=I_3(T^0_{\xi_3})^{+--}
$. The skew-symmetric part can be represented as $b_{\xi_i}=I_iu$, where $u$ is a traceless
symmetric (1,1)-tensor on $H$ which commutes with $I_1,I_2,I_3$. If $n=1$ then the tensor $u$
vanishes identically, $u=0$ and the torsion is a symmetric tensor, $T_{\xi}=T^0_{\xi}$.

As in \cite{IMV1},  we define two symmetric 2-tensors $T^0$ and $U$ on $H$ setting
\begin{equation}\label{tensors}
T^0(X,Y)=g((T_{\xi_1}^{0}I_1+T_{\xi_2}^{0}I_2+T_{%
\xi_3}^{0}I_3)X,Y),\quad U(X,Y)=g(u X,Y).
\end{equation}
It is easy to see that $T^0$ and $U$ are independent of the choice of the admissible set
$(\eta,\vartheta),$  and that they have the following properties:
\begin{equation}  \label{propt}
\begin{aligned} T^0(X,Y)+T^0(I_1X,I_1Y)+T^0(I_2X,I_2Y)+T^0(I_3X,I_3Y)=0, \quad tr_g(T^0)=tr_g(T^0I_s)=0 \\
U(X,Y)=U(I_1X,I_1Y)=U(I_2X,I_2Y)=U(I_3X,I_3Y),\quad tr_g(U)=tr_g(UI_s)=0. \end{aligned}
\end{equation}
In dimension seven $(n=1)$, the tensor $U$ vanishes identically,
$U=0$.

The  identity $ 4g(T^0(\xi_s,X),Y)=-T^0(I_sX,Y)-T^0(X,I_sY)$,  proved in
\cite[Proposition~2.3]{IV}, together with the first equality in \eqref{tensors} implies the
equivalence \cite{IMV1}
\begin{equation}\label{newequiv}
T^0=0 \Longleftrightarrow \{T^0_{\xi_s}=0, s=1,2,3\}.
\end{equation}
The torsion of Biquard connection is given in terms of the  tensors $T^0$ and $U$ by the formula
\begin{equation}\label{newtor}
g(T(\xi_s,X),Y)=g(T^0(\xi_s,X),Y)+U(I_sX,Y)=-\frac{T^0(I_sX,Y)+T^0(X,I_sY)}4+U(I_sX,Y).
\end{equation}
Let $R(X,Y)=[\nabla_X,\nabla_Y]-\nabla_{[X,Y]}$ be the curvature tensor of the Biquard connection.
The \emph{QC-Ricci curvature} $Ric$, the \emph{QC-Ricci forms} $\rho_s$ and the \emph{QC-scalar
curvature} $Scal$ are defined respectively by
$$Ric(A,B)=\sum_{a,b=1}^{4n}g(R(e_b,A)B,e_b), ~ A,B\in TM,$$
$$\rho_s(A,B)=\frac{1}{4n}\sum_{a=1}^{4n}g(R(A,B)e_a,I_s e_a),\quad Scal=\sum_{a,b=1}^{4n}
g(R(e_b,e_a)e_a,e_b),$$ where $e_1,...,e_{4n}$ is an orthonormal basis of $H$. The restriction of
the Ricci curvature $Ric$ to $H$ is a symetric 2-tensor (\cite{Biq1}) that could be
$Sp(n)Sp(1)$-invariantly decomposed in exactly three components. It is shown in (\cite{IMV1}) that
this three components are given by the 2-tensors $T^0,\ U$ and $Scal\cdot g.$ We have (see Theorem
3.12, \cite{IMV1}) :
\begin{equation} \label{Ric_components}
Ric(X,Y) \ =\ (2n+2)T^0(X,Y) +(4n+10)U(X,Y)+\frac{Scal}{4n}g(X,Y),\quad X,Y\in H.
\end{equation}

Since $V$ is preserved by $\nabla$ and $\nabla g=0$, there exist local 1-forms $\alpha_1$,
$\alpha_2$ and $\alpha_3$ such that
\begin{equation}\label{nabla-xi}
\nabla \xi_i=-\alpha_j\otimes \xi_k + \alpha_k\otimes \xi_j.
\end{equation}
Set $\tau=\displaystyle{\frac{Scal}{16n(n+2)}}$.
Then, according to \cite[Proposition 3.5~and~Theorem 3.12]{IMV1}, we have
\begin{equation}\label{alpha}
\alpha_i(\xi_s)=d\eta_s(\xi_j,\xi_k)-\delta_{is}(\tau+\frac{1}{2}d\eta_1(\xi_2,\xi_3)+\frac{1}{2}d\eta_2(\xi_3,\xi_1)
+\frac{1}{2}d\eta_3(\xi_1,\xi_2)).
\end{equation}

\section{The twistor space of a quaternionic contact manifold}

Let $(M,H,g,Q)$ be a quaternionic contact manifold. Let $\pi:Q\to M$ be the projection onto $M$ of
the bundle $Q$. Set
$$
{\mathcal Z}=\{I\in Q\, |\, I^2=-Id_H\}.
$$
Then $\pi_{\mathcal Z}=\pi|{\mathcal Z}:{\mathcal Z}\to M$ is a
subbundle of the vector bundle $Q$ called the {\it twistor space}
of the given QC-manifold. As we have mentioned, the condition $
I^2=-Id_H$ for $I\in Q$ is equivalent to $<I,I>=1$, thus
${\mathcal Z}=\{I\in Q\, |\,<I,I>=1\}$.

Let $\nabla$ be the Biquard connection on $M$ and denote by ${\mathcal H}$ the horizontal subbundle
of $TQ$ with respect to $\nabla$. For $I\in {\mathcal Z}$, the space ${\mathcal H}_I$ is tangent to
the submanifold ${\mathcal Z}$ of $Q$ since ${\mathcal Z}$ is the unit-sphere bundle of the vector
bundle $Q$ endowed with the metric $<.,.>$,  parallel with respect to the connection $\nabla$ on
$Q$. Further on, the restriction to ${\mathcal Z}$ of the horizontal bundle will be also denoted by
${\mathcal H}$. Let ${\mathcal V}$ be the vertical subbundle of $T{\mathcal Z}$. Then $T{\mathcal
Z}={\mathcal H}\oplus{\mathcal V}$.

 For $I\in {\mathcal Z}$, set $\xi_I=\varphi^{-1}(I)$ and denote by $\chi_I$ the horizontal lift of  $\xi_I$
at the point $I$. Fix a $I_0\in {\mathcal Z}$ and let
$(\eta,\vartheta)$ be an admissible set of the $QC$-structure
defined on a neighbourhood $U$ of the point $p=\pi(I_0)$. Denote
by $(\xi_1,\xi_2,\xi_3)$ the frame of vector fields on $U$
associated to $(\eta,\vartheta)$. Then every $I\in\pi_{\mathcal
Z}^{-1}(U)$ has a unique representation
$I=x_1(I)I_1+x_2(I)I_3+x_3(I)I_3$ where $x_1, x_2, x_3$ are smooth
functions such that $x_1^2+x_2^2+x_3^2=1$. We have
$\chi=x_1I^h_1+x_2I^h_3+x_3I ^h_3$ on $\pi_{\mathcal Z}^{-1}(U)$
where the upper script $h$ means "the horizontal lift". This shows
that $\chi$ is a smooth vector field on ${\mathcal Z}$.

 Let $g^h$ be the lift of the metric $g$ on $TM$ to the horizontal bundle ${\mathcal H}$ of ${\mathcal Z}$. Denote by
${\mathcal W}$ the orthogonal complement in ${\mathcal H}$ of the
horizontal vector field $\chi$. Then ${\mathcal D}={\mathcal
W}\oplus{\mathcal V}$ is a codimension $1$ subbundle of
$T{\mathcal Z}$ and, following \cite{Biq1}, we shall define an
almost complex structure ${\mathcal J}$ on it as follows.

Any vertical space ${\mathcal V}_I$, $I\in {\mathcal Z}$, is the tangent space at $I$ of the fibre
${\mathcal Z}_I$. The latter is the unit sphere in the $3$-dimensional vector space $Q_I$, so
${\mathcal V}_I=T{\mathcal Z}_I=\{S\in Q\,|\,<S,I>=0\}=\{S\in Q_I\,|\,SI+IS=0\}$. We define
${\mathcal J}|{\mathcal V}_I$ to be the standard complex structure of the $2$-sphere ${\mathcal
Z}_I$. In other words, we set ${\mathcal J}S=I\circ S$ for $S\in {\mathcal V}_I$.

 For $I\in {\mathcal Z}$, denote by $W_{I}$ the orthogonal complement  in $V_p$ of $\xi_I$, $V_p$ being the fibre
of $V$ at $p=\pi_{\mathcal Z}(I)$. We consider $W_{I}$ with the metric and the orientation induced
by those of $V_p$. Since the dimension of the space $W_{I}$ is $2$, there is a unique complex
structure $\hat{I}$ on it compatible with the metric and the orientation. If we denote by $\times$
the vector-cross product of the oriented Euclidean 3-dimensional vector space $V_p$, then
$\hat{I}{\zeta}=\xi_I\times\zeta$ for $\zeta\in W_I$. Note that the isomorphism $\varphi:V\to Q$
sends $W_I$ onto ${\mathcal V}_I\subset Q$ and
$$
\varphi(\hat{I}\zeta)={\mathcal J}\varphi(\zeta), \quad \zeta\in W_I.
$$
Now we define a complex structure $J_I$ on the space $H_p\oplus W_I$, $p=\pi_{\mathcal Z}(I)$,
setting $J_I|H_p=I$ and $J_I|W_I=\hat{I}$. Then we define ${\mathcal J}|{\mathcal W}_I$ as the
horizontal lift of $J_I$, i.e. ${\mathcal J}|{\mathcal W}_I$ is the pull-back of $J_I$ under the
isomorphism $\pi_{{\mathcal Z}\,\ast}:{\mathcal W}_I\to H_p\oplus W_I$.

 In this way we obtain a $CR$ manifold $({\mathcal Z},{\mathcal D},{\mathcal J})$. Recall that a Cauchy-Riemann ($CR$)
structure (in wide sense) on a manifold $N$ is a pair $({\mathcal D},{\mathcal J})$ of a subbundle
${\mathcal D}$ of the tangent bundle $TN$ and an almost complex structure ${\mathcal J}$ of the
bundle ${\mathcal D}$. For any two sections $X,Y$ of ${\mathcal D}$, the value of $[X,Y] \> mod\>
{\mathcal D}$ at a point $p\in N$ depends only on the values of $X$ and $Y$ at $p$,  and so we have
a skew-symmetric bilinear form $\omega: {\mathcal D}\times {\mathcal D}\to TN/{\mathcal D}$ defined
by $\omega(X,Y)=[X,Y] \> mod\> {\mathcal D}$; this form is called the Levi form of the
$CR$-structure $({\mathcal D},{\mathcal J})$. If the Levi form is ${\mathcal J}$-invariant, we can
define the Nijenhuis tensor of the $CR$-structure $({\mathcal D},{\mathcal J})$ by
$$
N^{\it CR}(X,Y)=-[X,Y]+[{\mathcal J}X,{\mathcal J}Y]-{\mathcal J}([{\mathcal J}X,Y]+[X,{\mathcal
J}Y]);
$$
its value at a point $p\in N$ lies in ${\mathcal D}$ and depends only on the values of the sections
$X,Y$ at $p$. A $CR$-structure is said to be integrable if its Levi form is ${\mathcal
J}$-invariant and the Nijenhuis tensor vanishes. Let ${\mathcal D}^{\mathbb C}={\mathcal
D}^{1,0}\oplus {\mathcal D}^{0,1}$ be the decomposition of the complexification of ${\mathcal D}$
into $(1,0)$ and $(0,1)$ parts with respect to ${\mathcal J}$. If the $CR$-structure $({\mathcal
D},{\mathcal J})$ is integrable, then the bundle ${\mathcal D}^{1,0}$ satisfies the following two
conditions:
$$
{\mathcal D}^{1,0}\cap \overline{{\mathcal D}^{1,0}}=0, \> [\Gamma({\mathcal
D}^{1,0}),\Gamma({\mathcal D}^{1,0})]\subset\Gamma({\mathcal D}^{1,0})
$$
where $\Gamma({\mathcal D}^{1,0})$ stands for the space of smooth
sections of ${\mathcal D}^{1,0}$. Conversely, suppose we are given
a complex subbundle ${\mathcal E}$ of the complexified tangent
bundle $T^{\mathbb C}N$ such that ${\mathcal E}\cap\overline
{\mathcal E}=0$ and $[\Gamma({\mathcal E}),\Gamma({\mathcal
E})]\subset\Gamma({\mathcal E})$ (many authors called a bundle
with these properties a "$CR$-structure"). Set ${\mathcal
D}=\{X\in TN: X=Z+\bar Z ~ \text{for some (unique)} ~ Z\in
{\mathcal E}\}$ and put ${\mathcal J}X=-{\it Im}Z$ for $X\in
{\mathcal D}$. Then $({\mathcal D},{\mathcal J})$ is an integrable
$CR$-structure such that ${\mathcal D}^{1,0}={\mathcal E}$.

 It a result of Biquard  \cite[Theorem II.5.1]{Biq1} that the $CR$-structure $({\mathcal D},{\mathcal J})$
on the twistor space ${\mathcal Z}$ is integrable and (up to isomorphism) is invariant under
conformal changes of the metric $g$.

Let $(\eta,\vartheta)$ be an admissible set of the given $QC$-structure defined on an open subset
$U$ of $M$. For $I=x_1I_1+x_2I_2+x_3I_3\in Z$, we set
$$
\eta^{\mathcal Z}_I=x_1\pi_{\mathcal Z}^{\ast}\eta_1 +x_2\pi_{\mathcal Z}^{\ast}\eta_2
+x_3\pi_{\mathcal Z}^{\ast}\eta_3.
$$
The right-hand side of the latter formula does not depend on
the choice of the admissible set $(\eta,\vartheta)$, thus we have a well-defined $1$-form
$\eta^{\mathcal Z}$ on ${\mathcal Z}$. It is clear that $\eta^{\mathcal Z}$ vanishes on ${\mathcal
V}$ and on the horizontal lift $H^h$ of the space $H$. Let $\xi_1,\xi_2,\xi_3$ be the frame of $V$
associated to $(\eta,\vartheta)$. Then $\xi_I=x_1\xi_1+x_2\xi_2+x_3\xi_3$, hence
$$\eta^{\mathcal Z}(\chi_I)=x_1^2+x_2^2+x_3^2=1.$$ Moreover, if $\zeta=\sum_{s=1}^3 z_s\xi_s\in W_I$, we have
$\eta^{\mathcal Z}(\zeta^h_I)=\sum_{s=1}^3 x_sz_s=0$. Therefore $Ker\,\eta^{\mathcal Z}={\mathcal
D}$. Define an endomorphism of the tangent bundle $T{\mathcal Z}$ setting $\Phi|{\mathcal
D}={\mathcal J}$ and $\Phi(\chi)=0$. Then
$$
\Phi^2(A)=-A+\eta^Z(A)\chi,\quad A\in T{\mathcal Z}.
$$
Thus $(\Phi,\chi,\eta^{\mathcal Z})$ is an almost contact structure on the twistor space ${\mathcal
Z}$ with contact distribution ${\mathcal D}$ and Reeb vector field $\chi$.

Let $(\eta,\vartheta)$ be an admissible set on an open subset $U$ of $M$ which is the domain of
local coordinates $u_1,...,u_m$ of $M$, $m=4n+3$. Set
$$
\widetilde u_r=u_r\circ\pi (I), r=1,...,m, \quad x_s(I)=<I,I_s>, s=1,2,3,
$$
for $I\in \pi^{-1}(U)\subset Q$. Then $(\widetilde u_1,...,\widetilde u_m,x_1,x_2,x_3)$ is a local
coordinate system of the manifold $Q$. For each vector field
$$X=\sum_{r=1}^m X^{r}\frac{\partial}{\partial u_r}$$
on $U$, the horizontal lift $X^h$ on $\pi^{-1}(U)$ is given by
\begin{equation}\label{hlift}
X^{h}=\sum_{r=1}^m (X^{r}\circ\pi)\frac{\partial}{\partial\tilde u_r}- \sum_{s,t=1}^3
x_s(<\nabla_{X}I_s,I_t>\circ\pi)\frac{\partial}{\partial x_t}.
\end{equation}


It follows from (\ref{hlift}) that
$$
[X^h,Y^h]=[X,Y]^h -\sum_{s,t=1}^3 x_s g(R(X,Y)I_s,I_t)\frac{\partial}{\partial x_t},
$$
where $R(X,Y)$ is the curvature tensor of  connection on $End(H)$ induced by the Biquard
connection. Let $I\in Q$ and $p=\pi(I)$. Using the standard identification of the tangent space
$T_IQ_I$ with the vector space $Q_I$ (the fibre of $Q$ through $I$), the latter formula can be
written as
\begin{equation}\label{bracket}
[X^h,Y^h]_I=[X,Y]_I^h -R_p(X,Y)I.
\end{equation}


\noindent {\bf Notation}. Fix a point $I\in{\mathcal Z}$ and set $p=\pi_{\mathcal Z}(I)$. Let
$I_1,I_2,I_3$ be an oriented orthonormal frame of $Q$ near the point $p$ such that $I_1(p)=I$ and
$\nabla I_s|_p=0$, $s=1,2,3$. If $(\eta',\vartheta')$ is any admissible set for the $QC$-structure
on $M$ near $p$, we have $I_s=\sum_{t=1}^3 a_{st}I_t'$ where $[a_{st}]\in SO(3)$. Then
$\eta=(\eta_1,\eta_2,\eta_3)$ with $\eta_s=\sum_{t=1}^3 a_{st}\eta_t'$ and
$\vartheta=(I_1,I_2,I_3)$ constitute an admissible set in a coordinate neighbourhood $U$ of $p$.
Let $(\widetilde u_1,...,\widetilde u_m,x_1,x_2,x_3)$ be the local coordinates  of the manifold $Q$
defined by means of this admissible set. Denote by $\xi_1,\xi_2,\xi_3$ the oriented orthonormal
frame of $V$ associated to $(\eta,\vartheta)$. Then $\nabla\xi_s|_p=0$, $s=1,2,3$.

  Given a section $a$ of $Q$, we denote by $\widetilde a$ the (local) vertical vector field on $Q$
defined by $\widetilde a_J=a_{\pi(J)}-<a_{\pi(J)},J>J$. This vector field is tangent to ${\mathcal
Z}$, thus its restriction to ${\mathcal Z}$, denoted again by $\widetilde a$, is a vertical vector
field on the twistor space.

 We shall use the above notation throughout this section and the following ones without further referring to it.

\begin{lemma}\label{L-H}
For any $I\in{\mathcal Z}$, a vector field $X$ on $M$ and a section $a$ of $Q$ near the point
$p=\pi_{\mathcal Z}(I)$, we have:
$$
\begin{array}{l}
[X^h,\widetilde a]_I=(\widetilde{\nabla_X a})_I, \quad
[\chi,\widetilde a]_I=-(\varphi^{-1}(a))^h_I +(\widetilde{\nabla_{\xi_I} a})_I\\[6pt]
[\chi,X^h]_I=(\nabla_{\xi_I}X+T_p(X,\xi_I))^h_I+R_p(X,\xi_I)I\\[6pt]
\end{array}
$$
\end{lemma}
\begin{proof}
Let $a=\sum_{s=1}^3 a_sI_s$. Then
$$
\widetilde a=\sum_{s=1}^3 \widetilde a_s\frac{\partial}{\partial x_s},~~\mbox{where}~~\widetilde
a_s=a_s-(\sum_{t=1}^3a_tx_t)x_s.
$$
We have
$$
X^h_I=\sum_{r=1}^m X^{r}(p)\displaystyle{\frac{\partial}{\partial\tilde u_r}}(I),\quad
[X^h,\displaystyle{\frac{\partial}{\partial x_s}}]_I=0,\quad \nabla_{X_p}a=\sum_{t=1}^3
X_p(a_t)I_t(p)
$$
since $\nabla I_s|_p=0$, $s=1,2,3$. Now the identity $[X^h,\widetilde a]_I=(\widetilde{\nabla_X
a})_I$ follows easily from  (\ref{hlift}). The second formula stated in the lemma follows from the
first one taking into account that (locally) $\chi=\sum_{s=1}^3 x_s\xi_s^h$. Moreover, in view of
(\ref{bracket}), we have
$[\chi,X^h]_I=\sum_{s=1}^3(x_s(I)[\xi_s,X]^h_I-R_p(\xi_s,X)I)=(\nabla_{\xi_I}X+T(X_p,\xi_I))^h_I+R_p(X,\xi_I)I$.
\end{proof}

\begin{lemma}\label{d-eta}
Let $I\in{\mathcal Z}$, $X,Y\in T_{\pi(I)}M$ and $a,b\in{\mathcal V}_I$. Let $X_H,Y_H$ and
$X_V,Y_V$ be, respectively, the $H$- and the $V$-components of $X,Y$. Then
$$
d\eta^{\mathcal Z}(X^h_I,Y^h_I)=2g(IX_H,Y_H)-2\tau g(\xi_I\times X_V,Y_V), \quad d\eta^{\mathcal
Z}(X^h_I,a)=-g(X,\varphi^{-1}(a)), \quad d\eta^{\mathcal Z}(a,b)=0.
$$
\end{lemma}
\begin{proof}
First, suppose that $X,Y\in H_p$, $p=\pi(I)$. Since $\nabla$ preserves $H$, we can extend the
vectors $X,Y$ to sections $X,Y$ of $H$ defined in a neighbouhood of $p$ such that $\nabla
X|_p=\nabla Y|_p=0$.

We have $X^h_I(\eta^{\mathcal
Z}(Y^h))=X^h_I(g^h(Y^h,\chi))=\sum_{s=1}^3X^h_I(x_sg(Y,\xi_s)\circ\pi)=\sum_{s=1}^3 x_s
X_p(g(Y,\xi_s))=0$ since $X^h_I(x_s)=0$, $\nabla Y|_p=0$ and $\nabla\xi_s|_p=0$, $s=1,2,3$.
Similarly, $Y^h_I(\eta^{\mathcal Z}(X^h))=0$. Recall that the form $\eta^{\mathcal Z}$ vanishes on
the vertical vectors. Then, in view of (\ref{bracket}), $d\eta^{\mathcal
Z}(X^h_I,Y^h_I)=-\eta^{\mathcal
Z}([X,Y]^h_I)=-\eta_1([X,Y]_p)=d\eta_1(X_p,Y_p)=2g_p(I_1X,Y)=2g_p(IX,Y)$.

 For any index $s=1,2,3$, we have $(\xi_s)^h_I(\eta^{\mathcal Z}(X^h))=0$ as above and $X^h_I(\eta^{\mathcal
 Z}(\xi_s^h))=\sum_{t=1}^3X^h_I(x_t\delta_{st})=0$. It follows that
$d\eta^{\mathcal Z}(X^h_I,(\xi_s)^h_I)=-\eta_1([X,\xi_s]_p)=\eta_1(T_p(X,\xi_s))=0$, $s=1,2,3$,
since $T_p(X,\xi_s)\in H$ (property $(vii)$ of the Biquard connection). Therefore $d\eta^{\mathcal
Z}(X^h_I,(\xi)^h_I)=0$ for $X\in H_p$ and $\xi\in V_p$.

We have also that $d\eta^{\mathcal
Z}((\xi_s)^h_I,(\xi_t)^h_I)=-\eta_1([\xi_s,\xi_t])_p=d\eta_1(\xi_s,\xi_t)_p$, $s,t=1,2,3$.  If
$\alpha_i$ are the $1$-forms defined by (\ref{nabla-xi}), then
$\alpha_i(\xi_s)=g(\nabla_{\xi_s}\xi_j,\xi_k)$ where $(i,j,k)$ is a cyclic permutation of
$(1,2,3)$. Therefore $\alpha_i(\xi_s)=0$ at the point $p$ and identity (\ref{alpha}) implies that
$d\eta_1(\xi_1,\xi_2)=d\eta_1(\xi_3,\xi_1)=0$ at $p$. Identity (\ref{alpha}) gives also that
\begin{equation}\label{2d}
2d\eta_i(\xi_j,\xi_k)-2\tau- d\eta_1(\xi_2,\xi_3)-d\eta_2(\xi_3,\xi_1) -
d\eta_3(\xi_1,\xi_2)=0~\mbox{ at the point}~p.
\end{equation}
Adding the three identities corresponding to the cyclic
permutations of $(1,2,3)$, we get from (\ref{2d}) that, at $p$,
$d\eta_1(\xi_2,\xi_3)+d\eta_2(\xi_3,\xi_1) +
d\eta_3(\xi_1,\xi_2)=-6\tau$. The latter identity and identity
(\ref{2d}) with $(i,j,k)=(1,2,3)$ give
$d\eta_1(\xi_2,\xi_3)=-2\tau$ at $p$. It follows that, for every
$\xi,\zeta\in V_p$, we have $d\eta^{\mathcal
Z}(\xi^h_I,\zeta^h_I)=-2\tau g(\xi_I\times\xi,\zeta)$.

Now let $X,Y\in T_pM$ be arbitrary tangent vectors. Writing $X=X_H+X_V$, $Y=Y_H+Y_V$ and applying
the preceding considerations we get the first formula stated in the lemma.

Next, take two sections of $Q$ with values $a$ and, respectively, $b$ at the point $p$, and zero
covariant derivatives at $p$. Denote these sections again by $a$ and $b$. Extend $X$ to a vector
field for which $\nabla X|_p=0$. Then  Lemma~\ref{L-H} implies that $d\eta^{\mathcal
Z}(X^h_I,a)=d\eta^{\mathcal Z}(X^h_I,\widetilde a_I)=-a(\eta^{\mathcal Z}(X^h))= -\sum_{s=1}^3
a(x_sg(X,\xi_s)\circ\pi)=-\sum_{s=1}^3 a(x_s)g_p(X,\xi_s)=-g(X,\varphi^{-1}(a)).$

  The last formula stated in the lemma follows from the fact that $\eta^{\mathcal Z}|{\mathcal
  V}=0$ and the bundle ${\mathcal V}$ is closed under the Lie bracket.
\end{proof}

\begin{cor}\label{d-A-chi}
$d\eta^{\mathcal Z}(A,\chi)=0$ for every $A\in T{\mathcal Z}$.
\end{cor}
\begin{proof}
If $A=X^h_I$ for a vector $X\in T_{\pi(I)}M$, then, by Lemma~\ref{d-eta}, $d\eta^{\mathcal
Z}(A,\chi)=-2\tau g(X_V\times\xi_I,\xi_I)=0$. If $A$ is a vertical vector, then $d\eta^{\mathcal
Z}(A,\chi)=g(\xi_I,\varphi^{-1}(A))=<\varphi(\xi_I),A>=<I,A>=0$.
\end{proof}

Set
\begin{equation}\label{defG}
G(A,B)=\frac{1}{2}d\eta^{\mathcal Z}(A,\Phi B)+\eta^{\mathcal Z}(A)\eta^{\mathcal Z}(B),\quad
A,B\in T{\mathcal Z}.
\end{equation}
Obviously, $G(\chi,\chi)=1$. By Corollary~\ref{d-A-chi}, $G(A,\chi)=G(\chi,A)=0$ for every $A\in
T{\mathcal Z}$. Moreover, it is easy to check by means of Lemma~\ref{d-eta} that $G$ is a symmetric
non-degenerate tensor. It is an observation of Biquard \cite[Theorem II.5.1]{Biq1} that the
symmetric form $d\eta^{\mathcal Z}(A,\Phi B)$ on the space $Ker\,\eta^Z$  is of signature
$(4n+2,2)$. Therefore $G$ is of signature $(4n+3,2)$; see also Corollary~\ref{G} below. Thus, $G$
is a pseudo-Riemannian metric on ${\mathcal Z}$ such that $G(\Phi A,\Phi B)=G(A,B)-\eta^{\mathcal
Z}(A)\eta^{\mathcal Z}(B)$ and $d\eta^{\mathcal Z}(A,B)=2G(\Phi A,B)$, $A,B\in T{\mathcal Z}$.
Therefore $(\Phi,\chi,\eta^{\mathcal Z},G)$ is a contact metric structure on the twistor space
${\mathcal Z}$.

\bigskip

\section{Normality of the contact structure on the twistor space of a $QC$-manifold}

The condition for the contact structure  $(\Phi,\chi,\eta^Z)$ to be normal is the vanishing of the
following tensor (cf., for example, \cite{Bl})

\begin{equation*}
N^1(A,B)=\Phi^2[A,B]+[\Phi A,\Phi B]-\Phi([\Phi A,B)]+[A,\Phi B])+2d\eta^Z(A,B)\chi,\quad A,B\in TZ.
\end{equation*}

Let $A,B$ be vector fields with values in $Ker\,\eta^Z={\mathcal D}$. Then
$$
N^1(A,B)=N^{\it CR}(A,B)+\eta^Z([A,B])\chi+d\eta^Z(A,B)\chi=N^{\it CR}(A,B).
$$
But, as we have mentioned, $N^{\it CR}(A,B)=0$ by a result of Biquard. Thus $N^1(A,B)=0$ for
$A,B\in Ker\,\eta^Z$ and to prove that $N^1=0$ it remains to show that $N^1(A,\chi)=0$ for $A\in
Ker\,\eta^Z$. The latter identity is equivalent to
$$
({\mathcal L}_{\chi}G)(A,B)=0 ~~ \mbox{for} ~~ A,B\in Ker\,\eta^Z={\mathcal D},
$$
where ${\mathcal L}$ stands for the Lie derivative. Indeed, we have  ${\mathcal
L}_{\chi}d\eta=(d\,\imath_{\chi}+\imath_{\chi}d)d\eta=d(\imath_{\chi}d\eta)=0$ since
$\imath_{\chi}d\eta=0$. On the other hand,  $({\mathcal L}_{\chi}d\eta)(A,\Phi B)=-\chi
G(A,B)+G([\chi,A],B)-G(A,\Phi[\chi,\Phi B])$ for $A,B\in {\mathcal D}$. Hence $\chi
G(A,B)-G([\chi,A],B)=-G(A,\Phi[\chi,\Phi B])$. Then $({\mathcal
L}_{\chi}G)(A,B)=-G(A,\Phi[\chi,\Phi B])-G(A,[\chi,B])=G(({\mathcal L}_{\chi}\Phi)(\Phi
B),A)=-G(N^1(B,\chi),A)$. This proves our claim since the form $G$ is non-degenerate.


 Lemma~\ref{d-eta} implies the following.
\begin{cor}\label{G}
Let $I\in{\mathcal Z}$, $X,Y\in T_{\pi(I)}M$ and $a,b\in{\mathcal V}_I$. Let $X_H,Y_H$ and
$X_V,Y_V$ be, respectively, the $H$- and the $V$-components of $X,Y$. Then
$$
\begin{array}{c}
G(X^h_I,Y^h_I)=g(X_H,Y_H)-\tau g(X_V,Y_V)+(\tau +1)g(X,\xi_I)g(Y,\xi_I),\\[6pt]
G(X^h_I,a)=\frac{1}{2}g(\xi_I\times X_V,\varphi^{-1}(a)),\quad G(a,b)=0.
\end{array}
$$
\end{cor}
\begin{proof}
We have
$Y^h_I=(Y_H)^h_I+(Y_V-g(Y,\xi_I)\xi_I)^h_I+g(Y,\xi_I)\chi$, hence
$\Phi Y^h_I=(IY)^h_I +(\xi_I\times Y_V)^h_I$ and $\eta^{\mathcal
Z}(Y^h_I)=g(Y,\xi_I)$. Now the statement follows from (\ref{defG})
and Lemma~\ref{d-eta}.
\end{proof}

  Corollary~\ref{G} can be also stated in the following form.
\begin{cor} Let $A,B\in T_{I_1}{\mathcal Z}$ and let
$$
A=X^h+\sum_{s=1}^3 u_s(\xi_s)^h+x_2I_2+x_3I_3,\quad B=Y^h+\sum_{s=1}^3 v_s(\xi_s)^h+y_2I_2+y_3I_3,
$$
where $X,Y\in H$. Then
$$
G(A,B)=g(X,Y)+u_1v_1-\tau(u_2v_2+u_3v_3)-\frac{1}{2}(u_3y_2+v_3x_2)+\frac{1}{2}(u_2y_3+v_2x_3).
$$
In particular, $G$ is of signature $(4n+2,2)$.
\end{cor}

We shall use Lemma~\ref{L-H} and Corollary~\ref{G} to prove the following.

\begin{lemma}\label{L-d-G}
Let $I\in{\mathcal Z}$, $X,Y\in H_p$, $p=\pi_(I)$, and $a,b\in{\mathcal V}_I$. Then
$$
\begin{array}{c}
({\mathcal L}_{\chi}G)(X^h_I,Y^h_I)=2g_p(T_{\xi_1}^0(X),Y),\quad ({\mathcal
L}_{\chi}G)(X^h_I,(\xi_s)^h_I)=\rho_s(X,\xi_1)_p+g_p([\xi_s,\xi_1],X),~ s=2,3,\\[6pt]
({\mathcal L}_{\chi}G)((\xi_s)^h_I,(\xi_t)^h_I)=-d\tau(\xi_1)_p\delta_{st}+\rho_s(\xi_t,\xi_1)_p+\rho_t(\xi_s,\xi_1)_p,~ s,t=2,3,\\[6pt]
({\mathcal L}_{\chi}G)(X^h,a)=0, \quad ({\mathcal L}_{\chi}G)((\xi_s)_I^h,a), ~ s=2,3,\\[6pt]
({\mathcal L}_{\chi}G)(a,b)=0.
\end{array}
$$
\end{lemma}

\begin{proof}
Extend the vectors $X, Y$ to sections $X, Y$ of $H$ for which $\nabla X|_p=\nabla Y|_p=0$.  By
Lemma~\ref{L-H} and Corollary ~\ref{G} we have
$$
({\mathcal
L}_{\chi}G)(X^h_I,Y^h_I)=\xi_1(g(X,Y))-g(T(X,\xi_1),Y)-g(T(Y,\xi_1),X)=2g(T_{\xi_1}^0(X),Y).
$$
Lemma~\ref{L-H} implies also that
$$
\begin{array}{l}
({\mathcal L}_{\chi}G)(X^h_I,(\xi_s)^h_I)=\xi^h_I(G(X^h,(\xi_s)^h))-G([\chi,X^h]_I,(\xi_s)^h_I)-G(X^h_I,[\chi,\xi_s^h]_I)\\[6pt]
=\xi^h_I(G(X^h,\xi_s^h))-G((T(X,\xi_I))^h_I,(\xi_s)^h_I)-G(R(X,\xi_I)I,(\xi_s)^h_I)\\[6pt]
-G((T(\xi_s,\xi_I))^h_I,X^h_I)-G(R(\xi_s,\xi_I)I,X^h_I).
\end{array}
$$
The first, second and last terms in the latter identity vanish by Corollary ~\ref{G} since $X\in H$
and $T(\xi_I,X)\in H$.  Moreover, the fourth term is equal to
$-g(T(\xi_s,\xi_1),X)=g([\xi_s,\xi_1],X)$ and for the third term we have, that when $s\neq 1$,
$$
\begin{array}{l}
G(R(X,\xi_I)I,(\xi_s)^h_I)=\displaystyle{\frac{1}{2}g_p(\xi_1\times\xi_s,\varphi^{-1}(R(X,\xi_1)I_1))
=\frac{1}{2}<R(X,\xi_1)I_1,I_1I_s>_p}=\\[6pt]
\displaystyle{\frac{1}{2}.\frac{1}{4n}(\sum_{a=1}^{4n}(-g(R(X,\xi_1)I_1e_a,I_sI_1e_a)
-g(I_1R(X,\xi_1)e_a,I_1I_se_a))}=-\rho_s(X,\xi_1),
\end{array}
$$
where $e_1,...e_{4n}$ is an orthonormal basis of $H_p$. This proves the second formula stated in
the lemma.

To prove the third formula, we note first that
$$
\begin{array}{c}
({\mathcal L}_{\chi}G)((\xi_s)^h_I,(\xi_t)^h_I)=\xi_1(-\tau\delta_{st})+\tau
g_p(T(\xi_s,\xi_1),\xi_t)-\frac{1}{2}g_p(\xi_1\times\xi_t,\varphi^{-1}(R(\xi_s,\xi_1)I_1))+\\[6pt]
g_p(T(\xi_t,\xi_1),\xi_s)-\frac{1}{2}g_p(\xi_1\times\xi_s,\varphi^{-1}(R(\xi_t,\xi_1)I_1))
\end{array}
$$
Next, $g_p(T(\xi_s,\xi_1),\xi_t)=-g_p([\xi_s,\xi_1],\xi_t)=d\eta_t(\xi_s,\xi_1)$ and, similarly,
$g_p(T(\xi_t,\xi_1),\xi_s)=d\eta_s(\xi_t,\xi_1)$. By (\ref{nabla-xi}) and (\ref{alpha}) we have

\begin{equation}\label{eto}
0=g_p(\nabla_{\xi_2}\xi_1,\xi_2)=\alpha_3(\xi_2)_p=d\eta_2(\xi_1,\xi_2)_p,\quad
0=g_p(\nabla_{\xi_3}\xi_3,\xi_1)=\alpha_2(\xi_3)_p=d\eta_3(\xi_3,\xi_1)_p\\[6pt]
\end{equation}
$$
0=g_p(\nabla_{\xi_2}\xi_3,\xi_1)+g_p(\nabla_{\xi_3}\xi_1,\xi_2)=\alpha_2(\xi_3)_p+\alpha_3(\xi_3)_p
=d\eta_2(\xi_3,\xi_1)_p-d\eta_3(\xi_1,\xi_2)_p.
$$
It follows that
$$
g_p(T(\xi_s,\xi_1),\xi_t)+g_p(T(\xi_t,\xi_1),\xi_s)=d\eta_t(\xi_s,\xi_1)_p+d\eta_s(\xi_t,\xi_1)_p=0,
\quad s,t=2,3.
$$
On the other hand
$$
\begin{array}{c}
\frac{1}{2}g_p(\xi_1\times\xi_t,\varphi^{-1}(R(\xi_s,\xi_1)I_1))=\frac{1}{2}<R(\xi_s,\xi_1)I_1,I_1I_t>_p
=-\rho_t(\xi_s,\xi_1)_p\\[6pt]
\frac{1}{2}g_p(\xi_1\times\xi_s,\varphi^{-1}(R(\xi_t,\xi_1)I_1))=\frac{1}{2}<R(\xi_t,\xi_1)I_1,I_1I_s>_p
=-\rho_j(\xi_t,\xi_1)_p.
\end{array}
$$
Thus
$$
({\mathcal
L}_{\chi}G)((\xi_s)^h_I,(\xi_t)^h_I)=-d\tau(\xi_1)\delta_{st}+\rho_s(\xi_t,\xi_1)+\rho_t(\xi_s,\xi_1).
$$

Now, extend the vertical vectors $a$ and $b$ to sections of the bundle $Q$. Then, by
Lemma~\ref{L-H} and Corollary ~\ref{G}, we have
$$
({\mathcal L}_{\chi}G)(X^h_I,a)=({\mathcal L}_{\chi}G)(X^h_I,\widetilde
a_I)=-G((T(X,\xi_I))^h_I,a)=0
$$
since $T(X,\xi_I)\in H_p$.

Applying Lemma~\ref{L-H} and Corollary ~\ref{G}, we get easily that
$$
({\mathcal L}_{\chi}G)((\xi_s)^h_I,a)=\xi_1(g(\xi_1\times \xi_s,\varphi^{-1}(\widetilde
a)))-\frac{1}{2}g_p(\xi_1\times
T(\xi_s,\xi_1)_V,\varphi^{-1}(a))-\tau(p)g_p(\xi_s,\varphi^{-1}(a)),\quad s=2,3.
$$
In the case when $a=I_t(p)$, $t=2,3$, the latter formula takes the form
$$
\begin{array}{c}
({\mathcal L}_{\chi}G)((\xi_s)^h_I,I_t(p))=\xi_1(g(\xi_1\times
\xi_s,\xi_t)+\frac{1}{2}g_p(\xi_1\times\xi_t,T(\xi_s,\xi_1))-\tau(p)g_p(\xi_s,\xi_t)=\\[6pt]
-\frac{1}{2}g_p(\xi_1\times\xi_t,[\xi_s,\xi_1])-\tau(p)\delta_{st}.
\end{array}
$$
If $s=t=2$ or $s=t=3$, then $g_p(\xi_1\times\xi_t,[\xi_s,\xi_1])$  is equal to
$d\eta_3(\xi_1,\xi_2)_p$ or to $d\eta_2(\xi_3,\xi_1)_p$, respectively. We have seen in the proof of
Lemma~\ref{d-eta} that $d\eta_1(\xi_2,\xi_3)=-2\tau$ at the point $p$. Similar simple arguments
show that $d\eta_2(\xi_3,\xi_1)=-2\tau$ and $d\eta_3(\xi_1,\xi_2)=-2\tau$ at $p$. Therefore
$({\mathcal L}_{\chi}G)((\xi_s)^h_I,I_t(p))=0$ for $s=t=2$ or $s=t=3$. This identity holds also for
$s=2, t=3$ and $s=3,t=2$ since $g_p(\xi_1\times\xi_t,[\xi_s,\xi_1])$ is equal to
$d\eta_2(\xi_2,\xi_1)_p=0$ by (\ref{eto}) in the first case and is equal to
$d\eta_3(\xi_1,\xi_3)_p=0$ in the second one. It follows that $({\mathcal
L}_{\chi}G)((\xi_s)^h_I,a)=0$.

Finally,
$$
\begin{array}{c}
({\mathcal L}_{\chi}G)(a,b)=\chi (G(\widetilde a,\widetilde b))-G([\chi,\widetilde
a)]_I,b)-G(a,[\chi,\widetilde b]_I)=\\[6pt]
\frac{1}{2}[g\xi_I\times\varphi^{-1}(a),\varphi^{-1}(b))+g(\xi_I\times\varphi^{-1}(b),\varphi^{-1}(a))=0.
\end{array}
$$
\end{proof}

\begin{thrm}\label{main_thrm}
The contact structure $(\Phi,\chi,\eta^{\mathcal Z})$ on the twistor space ${\mathcal Z}$ is ormal
if and only if the tensor $T^0$ on $M$ vanishes.
\end{thrm}

\begin{proof}
Recall that the normality condition for the structure $(\Phi,\chi,\eta^{\mathcal Z})$ is equivalent
to the condition $({\mathcal L}_{\chi}G)(A,B)=0 $ for $A,B\in Ker\,\eta^Z={\mathcal D}$. Thus, if
this structure is normal, then, by Lemma~\ref{L-d-G}, $T^{0}_{\xi_s}=0$, $s=1,2,3$. Hence, $T^0=0$
because of \eqref{newequiv}.

 Conversely, suppose that $T^{0}$=0. In view of Lemma~\ref{L-d-G}, to show that the structure $(\Phi,\chi,\eta^{\mathcal
 Z})$ is normal, we should prove that the following identities hold on $M$.
\begin{align}\label{MTE1}
&T^0_{\xi_1}=0; \\\label {MTE2} &\rho_{s}(X,\xi_1)_p+g_p([\xi_s,\xi_1],X)=0, \quad s=2,3,\quad X\in
H;\\\label {MTE3} &2\rho_s(\xi_s,\xi_1)_p=d\tau(\xi_1)_p,  \quad s=2,3;\\\label {MTE4}
&\rho_2(\xi_3,\xi_1)_p+ \rho_3(\xi_2,\xi_1)_p=0.
\end{align}
We shall prove that this system of equations follows from the single equation $T^0=0$.
The equation \eqref{MTE1} follows from \eqref{newequiv}.

To prove \eqref{MTE2}, note first that, according to  \cite[formula (4.3)]{IMV1}, we have
\begin{gather*}
\rho_i(X,\xi_j)=\frac{1}{2}d\eta_j([\xi_j,\xi_k]_H,X)=g(I_j[\xi_j,\xi_k]_H,X),\\
\rho_i(X,\xi_k)=\frac{1}{2}d\eta_k([\xi_j,\xi_k]_H,X)=g(I_k[\xi_j,\xi_k]_H,X),
\end{gather*}
where the subscript
$H$ means "the component in $H$". It follows that condition (\ref{MTE2}) is equivalent to the
identities
\begin{equation}\label{M}
I_1[\xi_3,\xi_1]_H=[\xi_1,\xi_2]_H,\quad I_1[\xi_1,\xi_2]_H=[\xi_1,\xi_3]_H
\end{equation}
It has been shown in \cite[formula (3.7)]{IV} that
\begin{equation}\label{vert023}
\begin{aligned}
3(2n+1)\rho_i(I_kX,\xi_j)=-\frac{(2n+1)(2n-1)}{16n(n+2)}X(Scal)
\hspace{2.3cm}\\
+\frac14\sum_{a=1}^{4n}(\nabla_{e_a}T^0)[(4n+1)(e_a,X)+3(I_ie_a,I_iX)]
+2(n+1)\sum_{a=1}^{4n}(\nabla_{e_a}U)(X,e_a).
\end{aligned}
\end{equation}
Thus we get
$$
3(2n+1)g(I_i[\xi_j,\xi_k],X)=-\frac{(2n+1)(2n-1)}{16n(n+2)}X(Scal)+2(n+1)\sum_{a=1}^{4n}(\nabla_{e_a}U)(X,e_a).
$$
The right hand-side of this identity does not depend on the indices $i,j,k$, therefore
$I_1[\xi_2,\xi_3]_H=I_2[\xi_3,\xi_1]_H=I_3[\xi_1,\xi_2]_H$ for every $(\xi_1,\xi_2,\xi_3)$. Then
$I_1[\xi_3,\xi_1]_H=-I_3I_2[\xi_3,\xi_1]_H=-I_3^2[\xi_1,\xi_2]_H=[\xi_1,\xi_2]_H$. Moreover,
writing the identity $I_1[\xi_2,\xi_3]_H=I_3[\xi_1,\xi_2]_H$ for $(\xi_3,\xi_1,\xi_2)$, we have
$I_3[\xi_1,\xi_2]_H=I_2[\xi_3,\xi_1]_H$, hence
$I_1[\xi_1,\xi_2]_H=I_2I_3[\xi_1,\xi_2]_H=I_2^2[\xi_3,\xi_1]_H=[\xi_1,\xi_3]_H$. This proves
(\ref{M}).


 Recall that, by definition,
$\rho_s(\xi_i,\xi_j)=\frac{1}{4n}\sum_{a=1}^{4n}g(R(\xi_i,\xi_j)e_a,I_se_a)$. According
to \cite[formula (3.6)]{IV},  we have
\begin{equation}\label{vert2}
\begin{array}{c}
g(R(\xi_i,\xi_j)e_a,I_se_a)=(\nabla_{\xi_i}U)(I_je_a,I_se_a)-(\nabla_{\xi_j}U)(I_ie_a,I_se_a)\\[8pt]
 -\frac14(\nabla_{\xi_i}T^0)(I_je_a,I_se_a)-\frac14(\nabla_{\xi_i}T^0)(e_a,I_jI_se_a)
 +\frac14(\nabla_{\xi_j}T^0)(I_ie_a,I_se_a)\\[8pt]+\frac14(\nabla_{\xi_j}T^0)(e_a,I_iI_se_a)]
 -(\nabla_{e_a}\rho_k)(I_iI_se_a,\xi_i) -\frac{Scal}{8n(n+2)}g(T(\xi_k,e_a),I_se_a)\\[8pt]
 -\sum_{b=1}^{4n}g(T(\xi_j,e_a),e_b)g(T(\xi_i,e_b),I_se_a)+\sum_{b=1}^{4n}g(T(\xi_j,e_b),I_se_a)g(T(\xi_i,e_a),e_b).
\end{array}
\end{equation}
In view of \eqref{newtor} and $T^0=0$, we obtain from \eqref{vert2} that
\begin{equation}\label{vert3}
\begin{array}{c}
g(R(\xi_i,\xi_j)e_a,I_se_a)=(\nabla_{\xi_i}U)(I_je_a,I_se_a)-(\nabla_{\xi_j}U)(I_ie_a,I_se_a)\\[8pt]
  -(\nabla_{e_a}\rho_k)(I_iI_se_a,\xi_i) -\displaystyle{\frac{Scal}{8n(n+2)}}U(I_ke_a,I_se_a)\\[8pt]
 -\sum_{b=1}^{4n}U(I_je_a,e_b)U(I_ie_b,I_se_a)+\sum_{b=1}^{4n}U(I_je_b,I_se_a)U(I_ie_a,e_b).
\end{array}
\end{equation}
Using the second identity in \eqref{propt}, we calculate
\begin{equation}\label{us}
\begin{aligned}
\sum_{a,b=1}^{4n}U(I_je_a,e_b)U(I_ie_b,I_ke_a)=\sum_{a,b=1}^{4n}U(I_je_a,e_b)U(e_b,I_je_a)=||U||^2,\\
\sum_{a,b=1}^{4n}U(I_je_a,e_b)U(I_ie_b,I_je_a)
=\sum_{a,b=1}^{4n}U(e_a,e_b)U(I_ie_b,e_a)=0.
\end{aligned}
\end{equation}
Now, identities (\ref{vert3}), \eqref{us} and the fact that $U$ is completely trace-free imply that
\begin{align}
\label{rho-i}
&4n\rho_i(\xi_i,\xi_j)=\sum_{a=1}^{4n}g(R(\xi_i,\xi_j)e_a,I_ie_a)=\sum_{a=1}^{4n}(\nabla_{e_a}\rho_k)(e_a,\xi_i)
\\ \label{rho-j}
&4n\rho_j(\xi_i,\xi_j)=\sum_{a=1}^{4n}g(R(\xi_i,\xi_j)e_a,I_je_a)=-\sum_{a=1}^{4n}(\nabla_{e_a}\rho_k)(I_ke_a,\xi_i)\\
\label{rho-k}
&4n\rho_k(\xi_i,\xi_j)=\sum_{a=1}^{4n}g(R(\xi_i,\xi_j)e_a,I_ke_a)
=\sum_{a=1}^{4n}(\nabla_{e_a}\rho_k)(I_je_a,\xi_i)-2||U||^2.
\end{align}

Next, according to (\ref{vert023}), we have
\begin{equation}\label{rhoxi}
3(2n+1)\rho_k(X,\xi_i)=(2n+1)(2n-1)d\tau(I_jX)+2(n+1)\sum_{a=1}^{4n}(\nabla_{e_a}U)(I_jX,e_a),
\quad X\in H.
\end{equation}
Let $X\in H_p$. Extend $X$ to a local vector field such that $\nabla X|_p=0$. Then we get from
(\ref{rhoxi}) that
$$
3(2n+1)(\nabla_X\rho_k)(X,\xi_i)=(2n+1)(2n-1)X(d\tau(I_jX))+2(n+1)\sum_{a=1}^{4n}X((\nabla_{e_a}U)(I_jX,e_a)).
$$
This identity and identities (\ref{rho-i}), (\ref{rho-j}), (\ref{rho-k}) imply that at the point
$p$ we have:
\begin{multline}\label{rho-iij}
12(2n+1)\rho_i(\xi_i,\xi_j)=\\\sum_{a=1}^{4n}[(2n+1)(2n-1)e_a(d\tau(I_je_a))
+2(n+1)\sum_{a,b=1}^{4n}e_a((\nabla_{e_b}U)(I_je_a,e_b))],
\end{multline}
\begin{multline}\label{rho-jij}
12(2n+1)\rho_j(\xi_i,\xi_j)=\\-\sum_{a=1}^{4n}[(2n+1)(2n-1)e_a(d\tau(I_ie_a))
+2(n+1)\sum_{a,b=1}^{4n}e_a((\nabla_{e_b}U)(I_ie_a,e_b))]\end{multline}
\begin{multline}\label{rho-kij}
3(2n+1)4n\rho_k(\xi_i,\xi_j)=\\-\sum_{a=1}^{4n}[(2n+1)(2n-1)e_a(d\tau(e_a))
+2(n+1)\sum_{a,b=1}^{4n}e_a((\nabla_{e_b}U)(e_a,e_b))]-6(2n+1)||U||^2.
\end{multline}
The right hand-side of the last identity does not depend on the indices $i,j,k$, therefore
$\rho_3(\xi_1,\xi_2)_p=\rho_2(\xi_3,\xi_1)_p$, which proves (\ref{MTE4}).

 Identities (\ref{rho-iij}) and (\ref{rho-jij}) imply that
$\rho_3(\xi_3,\xi_1)_p=-\rho_2(\xi_1,\xi_2)_p$, and so
$\rho_3(\xi_3,\xi_1)_p=\rho_2(\xi_2,\xi_1)_p$. On the other hand, we have
$\rho_3(\xi_3,\xi_1)+\rho_2(\xi_2,\xi_1)=\xi_1(\tau)$ by \cite[formula (4.6)]{IMV1}. Therefore
$2\rho_2(\xi_2,\xi_1)_p=2\rho_3(\xi_3,\xi_1)_p=d\tau(\xi_1)_p$, i.e. condition (\ref{MTE3}) also
holds.
\end{proof}

Theorem~\ref{main_thrm}, (\ref{Ric_components}) and (\ref{propt}) imply the following

\begin{thrm}
The contact structure $(\Phi,\chi,\eta^{\mathcal Z})$ on the twistor space ${\mathcal Z}$ is normal
if and only if the QC Ricci tensor commutes with the quaternionic structure on the contact
distribution, $$Ric(IX,IY)=Ric(X,Y), \quad X,Y\in H, I\in Q.$$
\end{thrm}


\begin{thebibliography}{99}

\bibitem{AFIV} de Andres, L., Fernandez, M., Ivanov, S., Joseba, S. \&
Ugarte, L. \& Vassilev, D., {\em  Explicit  Quaternionic Contact
Structures and Metrics with Special Holonomy },  arXiv:0903.1398
v3.

\bibitem{AFIV1} de Andres, L., Fernandez, M., Ivanov, S., Joseba, S. \&
Ugarte, L. \& Vassilev, D., {\em  Quaternionic Kaehler and Spin(7) metrics
arising from quaternionic contact Einstein structures},   arXiv:1009.2745.

\bibitem{Biq1} Biquard, O., \emph{M\'{e}triques d'Einstein
asymptotiquement sym\'{e}triques.} Ast\'{e}risque \textbf{265} (2000), vi+109 pp.

\bibitem{Biq2} Biquard, O., \emph{Quaternionic contact structures.}
Quaternionic structures in mathematics and phisics (Rome, 1999), 23--30 (electronic), Univ. Studi
Roma "La Sapienza", Roma, 1999.

\bibitem{Bl} Blair, D., \emph{Reimannian geometry of contact and symplectic manifolds}, Progress in
Mathematics, vol. 203, Birkh\"auser, Boston-Basel-Berlin, 2002.


\bibitem{D} Duchemin, D.,  \emph{Quaternionic contact structures in
dimension 7.} Ann. Inst. Fourier (Grenoble) \textbf{56} (2006), no. 4, 851--885.

\bibitem{D1} Duchemin, D.,  \emph{Quaternionic contact hypersurfaces.}
math.DG/0604147.



\bibitem{F2} Folland, G.,  \emph{Subelliptic estimates and function spaces on nilpotent Lie groups},
  Ark. Math., \textbf{13}~(1975), 161--207.


\bibitem{FS} Folland, G., \& Stein, E., \ \emph{Estimates for the $\Bar {\partial}_{b}$ Complex and Analysis on
the Heisenberg Group}, Comm. Pure Appl. Math., \textbf{27}~(1974),
429--522.



\bibitem{GV}
Garofalo, N., \& Vassilev, D.,\ \emph{ Symmetry properties of
positive entire solutions of Yamabe type equations on groups of
Heisenberg type}, Duke Math J, {\bf 106} (2001), no. 3, 411--449.



\bibitem{IMV1} Ivanov, S., Minchev, I., \& Vassilev, D.,\ \emph{%
Quaternionic contact Einstein structures and the quaternionic contact Yamabe problem.} preprint,
math.DG/0611658.

\bibitem{IMV2} Ivanov, S., Minchev, I., \& Vassilev, D.,\ \emph{Extremals for the Sobolev
inequality on the seven dimensional quaternionic Heisenberg group and the quaternionic contact
Yamabe problem}, J. Eur. Math. Soc., {\bf 12} (2010), 1041–-1067.

\bibitem{IMV3} Ivanov, S., Minchev, I., \& Vassilev, D., \ \emph{work in progress}

\bibitem{IMV4} Ivanov, S., Minchev, I., \& Vassilev, D.,\ \emph{The optimal constant in the $L^2$ Folland-Stein
inequality on the quaternionic Heisenberg group },
math.DG/1009.2978.

\bibitem{IV}  Ivanov, S.,  \&  Vassilev, D.,\ \emph{Conformal quaternionic contact curvature and the local
sphere theorem},  J. Math. Pures Appl., {\bf 93} (2010), 277-307.


\bibitem{IV1} Ivanov, S.,  \& Vassilev, D.,\ \emph{Quaternionic contact manifolds with a closed fundamental
4-form},    Bull. London Math. Soc. \textbf{42} (2010), 1021-1030.

\bibitem{LeB} LeBrun, C.,  {\em $\mathcal H$-spaces with a cosmological
constant}, Proc. Roy.Soc. London Ser. A \textbf{380} (1982),
171-185

\bibitem{LeB2} LeBrun, C., {\em On complete quaternionic-K\"ahler manifolds},
Duke Math. J. \textbf{63} (1991), 723-743

\bibitem{Wei} Wang, W., \emph{The Yamabe problem on quaternionic contact
manifolds}, Ann. Mat. Pura Appl., \textbf{186} (2007), no. 2, 359--380.





\end{thebibliography}
\end{document}